# An Iteratively Reweighted Algorithm for Sparse Reconstruction of Subsurface Flow Properties from Nonlinear Dynamic Data

Lianlin Li and B. Jafarpour

---

A challenging problem in predicting fluid flow displacement patterns in subsurface environment is the identification of spatially variable flow-related rock properties such as permeability and porosity. Characterization of subsurface properties usually involves solving a highly underdetermined nonlinear inverse problem where a limited number of measurements are used to reconstruct a large number of unknown parameters. To alleviate the non-uniqueness of the solution, prior information is integrated into the solution. Regularization of ill-posed inverse problems is commonly performed by imparting structural prior assumptions, such as smoothness, on the solution. Since many geologic formations exhibit natural continuity/correlation at various scales, decorrelating their spatial description can lead to a more compact or sparse representation in an appropriate compression transform domain such as wavelets or Fourier domain. The sparsity of flow-related subsurface properties in such incoherent bases has inspired the development of regularization techniques that attempt to solve a better-posed inverse problem in these domains. In this paper, we present a practical algorithm based on sparsity regularization to effectively solve nonlinear dynamic inverse problems that are encountered in subsurface model calibration. We use an iteratively reweighted algorithm that is widely used to solve linear inverse problems with sparsity constraint (known as compressed sensing) to estimate permeability fields from nonlinear dynamic flow data. To this end, we minimize a data misfit cost function that is augmented with an additive regularization term promoting sparse solutions in a Fourier-related discrete cosine transform domain. This regularization approach introduces a new weighting parameter that is in general unknown a priori, but controls the effectiveness of the resulting solutions. Determination of the regularization parameter can only be achieved through considerable numerical experimentation and/or a *priori* knowledge of the reconstruction solution. To circumvent this problem, we include the sparsity promoting constraint as a multiplicative regularization term which eliminates the need for a regularization parameter. We evaluate the performance of the iteratively reweighted approach with multiplicative sparsity regularization using a set of waterflooding experiments in an oil reservoir where we use nonlinear dynamic flow data to infer the spatial distribution of rock permeabilities, While, the examples of this paper are derived from the subsurface flow and transport application, the proposed methodology also can be used in solving nonlinear inverse problems with sparsity constraints in other imaging applications such as geophysical, medical imaging, electromagnetic and acoustic inverse problems.

**Key words:**
History matching, oil exploration, the sparsity constraint, iterative reweighed approach, the artificial regularization parameters, multiplicative regularization, compressive sensing

# I.  INTRODUCTION

Mathematical models governing fluid flow and transport in porous media are widely used to quantify and predict fluid displacement behavior in the subsurface environment. These models require specification of subsurface fluid and rock properties as input parameters to derive the solution for the spatial and temporal distribution of the system's states (such as pressure and saturation distribution). In general, the reliability of model predictions depends on the assumptions used in generating the models and accuracy of the input parameters used to obtain the solution. Of particular importance in determining the behavior of fluid flow and displacement is the rock permeability distribution, which explains why characterization of rock permeability is an essential part of any development and management plans for exploitation of groundwater and hydrocarbon resources as well as site remediation and cleanup activities. Since permeability can not be sampled directly and because data collection is very costly, the spatial distribution of permeability is usually inferred from indirect measurements through a process known as inverse modeling (or parameter estimation, model calibration, history matching).

The inverse problem of estimating heterogeneous permeability distribution from nonlinear dynamic data (such as, pressure and flow measurements) that are usually taken at scattered injection/production wells is ill-posed. For heterogeneous reservoirs, the history matching solution is generally non-unique due to highly undetermined nature of the problem, i.e. there are significantly larger number of unknown parameters than the available data. As a consequence, several different solutions may be found that are consistent with all observed measurements but fail to forecast the behavior of the fluid displacement correctly.

The underlying inverse problem is usually formulated to minimize a prescribed cost function that describes the misfit between measured and simulated data by tuning the unknown parameter values (permeability in this case) while honoring the governing flow equations. Because of ill-posedness, a form of prior information has to be incorporated into the solution procedure to eliminate solutions that are inconsistent with the specified prior knowledge. Prior information can be incorporated in the solution either in the form of an explicit set of initial parameters that should not be deviated from or an implicit structural characteristic imposed on the solution (such as smoothness). In any case, while prior information is needed to constrain the solution cautions must be exercised against placing emphasis on uncertain prior knowledge. Another common approach for regularizing ill-posed inverse problem is achieved through suitable parameterization of the properties. When conducted correctly, parameterization can reduce of the number of model parameters with minimal effect on solution accuracy while simultaneously improving the computational complexity of the solution approach.

A proper regularization approach should be based on the physics of the problem and the available qualitative or quantitative knowledge about the expected features in the solution. In a geologic setting, for instance, a prograding formation with gradually varying heterogeneous property may warrant smoothness as a structural regularization constraint while in a channelized or fractured medium where sudden discontinuities are expected, smoothness can not be justified as a regularization constraint. Although spatially driven structural assumptions have been successfully used in several geophysical applications, the intrinsic attributes of geologic formations may be exploited to derive more effective structural priors in an appropriate transform domain. Therefore, searching for more effective regularization techniques for subsurface properties is appealing.

The intrinsic continuity (spatial correlation) in geologic formations are amenable to a high degree of decorrelation and energy compaction that can be achieved through the existing compression transforms such as the discrete cosine basis, wavelet basis, and many other prereconstructed or learned (trained) dictionaries and so on that allow for sparse representation of the most salient information in pixel-based spatial descriptions. For example, continuous geologic features, such as facies descriptions in the spatial domain translate into a sparse representation in an appropriate transform domain. Jafarpour et al. (2009a) exploited this property to formulate and apply a two-step sparsity regularization algorithm to reconstruct facies distribution from production data. In this paper, by taking advantage of the sparsity, the key concept in the compressive sensing theory, in a transform domain representation of permeability distribution, we develop a novel approach for reconstruction of the unknown high-dimensional permeability distribution from limited measurements. Our development in here has several important differences with those in [REF]. Unlike in [REF] we do not search for a solution in a truncated DCT subspace and include the complete DCT basis in the solution. We use a more efficient iteratively reweighted solution algorithm that imposes the sparsity constraint more effectively. Moreover, enforce the sparsity regularization through a multiplicative (rather than an additive) term and consider reconstruction under noisy measurements. Although the additive regularization term has been shown to improve the solution quality, it introduces an additional regularization parameter in the cost function that can only be determined through considerable numerical experimentation and/or using a priori knowledge of the reconstruction solution. We show how a multiplicative regularization term can be used to circumvent this issue. While a multiplicative term may increase the nonlinearity of the objective function, our results indicate that this approach effectively imposes the sparsity regularization without needing an additional regularization parameter. Based on the examples in this paper, the multiplicative regularization approach seems to reconstruct permeability fields from limited noisy data without any a priori information in addition to sparsity of the solution in the transform domain.

The reminder of this paper is organized as follows. Section II describes the formulation to solve the dynamic inverse problem with sparsity constraint. The iteratively reweighted algorithm with additive and multiplicative sparsity-promoting regularization is presented in Section III. The reconstruction results for several two-dimensional examples of waterflooding in an oilfield are provide in Section IV, followed by the discussion of the results and presentation of conclusions in Section V.

## II. PROBLEM STATEMENT

We begin the problem formulation by briefly introducing the notations and assumptions that we have used. Let $J_n(\mathbf{m})$ be the cost-function to be optimized at iteration $n$ with $\mathbf{m} \in \mathbb{R}^N$ being the vector of unknown parameters (permeabilities in each grid block) to be determined, and $\mathbf{m}^{(n)} \in \mathbb{R}^N$ indicating the unknown parameters at iteration $n$. We denote the vector of measured data (pressure, water saturation, etc.) with $\mathbf{y} \in \mathbb{R}^M$ and use $\mathbf{y}^{sim} = \mathbf{g}(\mathbf{m})$ to distinguish the vector of simulated data (given that the parameters $\boldsymbol{m}$) from the measurements. The simulated data on injection and/or production wells are nonlinear functions of the parameters $\boldsymbol{m}$; thus, making $J$ a nonlinear cost-function. Moreover, we use $\|x\|_p^p = \sum_{i=1}^{N} |x_i|^p$ ($0 < p < \infty$) to denote the $lp$-norm of vector $\boldsymbol{x}$, and use the classical iterative Newton approach to optimize the specified cost-function.

A common cost-function to minimize with respect to $\boldsymbol{m}$ is $J(\mathbf{m}) = \|\mathbf{y} - \mathbf{g}(\mathbf{m})\|_2^2$. In ill-posed inverse problems, the resulting linear equations in the Newton's method can be highly ill-conditioned during the iterations, which can be improved using regularization techniques. A common form of regularization uses the minimal energy solution or the well-known $l_2$-norm

constraint, i.e. $\|\mathbf{m}\|_2^2$. The resulting cost-function in this case is

$$J(\mathbf{m}) = \|\mathbf{y} - \mathbf{g}(\mathbf{m})\|_2^2 + \alpha \|\mathbf{m}\|_2^2 \tag{1}$$

where $\alpha$ is a regularization parameter that balance the data misfit and the regularization term. While the above minimum energy regularization have been successfully applied to many problems [REFS], it does not have sparsity-promoting property. To impose sparsity on the solution an alternative regularization term must be used. The relevance of sparsity as a regularization constraint that we exploit in this paper can be better appreciated by referring to Fig. 1. Fig. 1a shows a typical permeability distribution with its corresponding DCT coefficients in Fig. 1b. An approximate reconstruction of the permeability field in Fig. 1a using only 5% largest DCT coefficients is shown in Fig. 1c .

In an inverse problem where we do not know the true solution, the location and magnitude of the significant coefficients have to be inferred from available measurements. Since the problem is ill-posed, we exploit the sparsity of the DCT-domain representation of the property image to help find a solution. For linear problems, the sparse reconstruction techniques such as the compressed sensing paradigm outline the solution approach. However, the underlying concepts in the linear case can be used to formulate solution procedures for nonlinear problems with dynamic measurements.

The original sparse reconstruction problem is formulated as

$$, \min_{\mathbf{m}} \|\mathbf{\Phi}\mathbf{m}\|_0 \tag{P0}$$

$$s.t. \quad \mathbf{y} = \mathbf{g}(\mathbf{m}) \quad \text{or} \quad \|\mathbf{y} - \mathbf{g}(\mathbf{m})\|_2^2 \leq \sigma$$

where $\sigma$ is the noisy energy and the $l_0$-quasinorm minimization attempts to find solutions with minimum support that satisfy the measurements errors within the specified noise. The exact solution to (P0) requires exhaustive search over all possible sparse set of the basis and is combinatorial (i.e. NP-hard). However, there a vast literature approximate solution to (P0) and elegant theoretical results, under mild conditions, have been obtained for the case where $g(\boldsymbol{m})$ is linear. The practical solution techniques can be classified as methods that apply a *convex relaxation* [REF] to (P0) and solve the corresponding convex optimization problem and *greedy algorithms* known as matching pursuit [REF]. A generalized form of the original problem in (P0) can be formulated in two ways. The first formulation is written as a constrained nonlinear optimization problem (P1)

$$\min_{\mathbf{m}} \left\| \mathbf{\Phi m} \right\|_p^p \tag{P1}$$

$$s.t. \quad \mathbf{y} = \mathbf{g}(\mathbf{m}) \quad or \quad \left\| \mathbf{y} - \mathbf{g}(\mathbf{m}) \right\|_2^2 \leq \sigma \; .$$

.whereas the second formulation is stated as an unconstrained nonlinear optimization problem (P2) as follows

$$\min_{\mathbf{m}} \mathcal{L}(\mathbf{m}) = \left\| \mathbf{y} - \mathbf{g}(\mathbf{m}) \right\|_2^2 + \alpha \left\| \mathbf{\Phi m} \right\|_p^p \tag{P2}$$

where $\alpha$ is the regularization parameter balancing the data misfit and regularization term.

In this paper, we focus on the efficient algorithm to solve (P2) which will be treated with more details in the next section. The remainder of this section is focused on a brief discussion of (P1). The Lagranian $\mathcal{L} : \mathbb{R}^M \times \mathbb{R}^M \rightarrow \mathbb{R}$ associated with problem (P1) can be formulated as

$$L(\mathbf{m}, \boldsymbol{\lambda}) = \left\| \mathbf{\Phi m} \right\|_p^p + \boldsymbol{\lambda}^T \left( \mathbf{y} - \mathbf{g}(\mathbf{m}) \right) \tag{2}$$

where $\boldsymbol{\lambda} \in \mathbb{R}^M$ is the Lagrange multiplier associated with the equality constraint $\mathbf{y} = \mathbf{g}(\mathbf{m})$. The vector $\boldsymbol{\lambda} \in \mathbb{R}^M$ is called the dual variables or Lagrange multiplier vector associated with the problem (P1). The optimal solution $(\tilde{\mathbf{m}}, \tilde{\boldsymbol{\lambda}})$ to (2) can be obtained by solving the following normal equations

$$\nabla_m \mathcal{L}(\mathbf{m}, \boldsymbol{\lambda}) = p\mathbf{\Phi}^T \mathbf{D} \mathbf{\Phi m} + \mathbf{G}^T \boldsymbol{\lambda} = \mathbf{0} \tag{3}$$

$$\nabla_\lambda \mathcal{L}(\mathbf{m}, \boldsymbol{\lambda}) = \mathbf{y} - \mathbf{g}(\mathbf{m}) = \mathbf{0} \tag{4}$$

where $D$ is the diagonal matrix with entries $\mathbf{D}_{i,i} = (\mathbf{\Phi m})_i$ and $\mathbf{G}$ is the sensitivity matrix which can be efficiently computed using the adjoint method (REF). To solve the nonlinear equations in (3) and (4), linearized approximation with iterative solution schemes can be used. Using an iterative framework, at iteration $n+1$, $\mathbf{g}(\mathbf{m})$ is assumed to be approximated by $\mathbf{g}(\mathbf{m}^{(n)}) + \mathbf{G}(\mathbf{m} - \mathbf{m}^{(n)})$, and the entries of $\mathbf{D}$ are approximated by $\mathbf{D}_{i,i} = (\mathbf{\Phi m}^{(n)})_i$. From Equation (3), the solution is written as

$$\mathbf{m} = -\frac{1}{p} (\mathbf{\Phi}^T \mathbf{D} \mathbf{\Phi})^{-1} \mathbf{G}^T \boldsymbol{\lambda} \tag{5}$$

$$\boldsymbol{\lambda} = -p \left( \mathbf{G} \left( \boldsymbol{\Phi}^T \mathbf{D} \boldsymbol{\Phi} \right)^{-1} \mathbf{G}^T \right)^{-1} \mathbf{y}_n, \tag{6}$$

where $\mathbf{y}_n = \mathbf{y} - \mathbf{g}\left(\mathbf{m}^{(n)}\right) + \mathbf{G}\mathbf{m}^{(n)}$. After combining Equations (5) and (6), the solution at iteration $n$+1 is

$$\mathbf{m} = \left( \boldsymbol{\Phi}^T \mathbf{D} \boldsymbol{\Phi} \right)^{-1} \mathbf{G}^T \left( \mathbf{G} \left( \boldsymbol{\Phi}^T \mathbf{D} \boldsymbol{\Phi} \right)^{-1} \mathbf{G}^T \right)^{-1} \mathbf{y}_n \tag{7}$$

Solving (P1) using the presented formulation can be nontrivial because $\mathbf{G}\left( \boldsymbol{\Phi}^T \mathbf{D} \boldsymbol{\Phi} \right)^{-1} \mathbf{G}^T$ is in general highly ill-posed. A discussion of the solution approach for (P1) is beyond the scope of the current paper.

In this paper, we have derived the adjoint formulation for a two-phase reservoir, starting directly from the governing nonlinear partial differential equations. The adjoint model is solved through a backward recursion after solving the forward model as described in [REF]. In solving the inverse problem we start with a *homogeneous* permeability distribution and calculate the difference between the computed and the given production data while simoultaneously minimizing the sparsity regularization. These residuals are minimized by using the computed gradients from the adjoint model. The procedure is repeated iteratively until a stopping (convergence) criterion is satisfied.

### III. DESCRIPTION OF PROPOSED ALGORITHM

In this section we describe an efficient algorithm for solving (P2) via classical iterative Newton method. It has been shown [REF] that the quasinorm $\left\| \boldsymbol{\Phi}\mathbf{m} \right\|_p^p$ in (P2) imparts sparsity on the solution only for $0 < p \leq 1$. However, for this range of $p$ the gradient of $\left\| \boldsymbol{\Phi}\mathbf{m} \right\|_p^p$ with respect to $m$ is a highly non-smooth function of $m$ and is not defined at zero. Consequently, (P2) is a non-smooth optimization problem and sub-gradient based solution methods are not likely to perform well. In this paper we use the iteratively reweighted approach that is widely used in the field of compressive sensing [REF] to minimize (P2). As mentioned above, augmenting the misfit term with the sparsity constraint $\left\| \boldsymbol{\Phi}\mathbf{m} \right\|_p^p$ introduces an additional regularization parameter $\alpha$, which is in general not easy to specify without considerable numerical experimentation and *a priori* information about the solution. To eliminate the need to specify $\alpha$, we include the regularization term by multiplying it to the misfit function as described in Section III.B.

## III.A Additive Regularization

In this regularization approach, the cost-functional to be minimized at ($n$+1)-th iteration is

$$J_{n+1}(\mathbf{m}) = \left\| \mathbf{y} - \mathbf{g}(\mathbf{m}) \right\|_2^2 + \alpha \left\| \boldsymbol{\Phi}\mathbf{m} \right\|_{\mathbf{W}}^2 \tag{8}$$

where $\alpha$ is a regularization parameter that balances the weight given to the prior information $\left\| \boldsymbol{\Phi}\mathbf{m} \right\|_{\mathbf{W}}^2$ and the data misfit $\left\| \mathbf{y} - \mathbf{g}(\mathbf{m}) \right\|_2^2$. In our formulation, we use the the regularization term $\left\| \boldsymbol{\Phi}\mathbf{m} \right\|_{\mathbf{W}}^2 = \mathbf{m}^T \boldsymbol{\Phi}^T \mathbf{W} \boldsymbol{\Phi}\mathbf{m}$, where $\mathbf{W}$ is a diagonal "weighting" matrix with entries $\mathbf{W}_{i,i}$ approximately proportional to $\left(\boldsymbol{\Phi}\mathbf{m}^{(n)}\right)_i^{p-2}$, is used as the regularization term. To avoid the singularity due to very small values of $\left(\boldsymbol{\Phi}\mathbf{m}^{(n)}\right)_i$, $\mathbf{W}_{i,i}$ is modified into $\mathbf{W}_{i,i} = \left(\left(\boldsymbol{\Phi}\mathbf{m}^{(n)}\right)_i^2 + \varepsilon_n\right)^{\frac{p}{2}-1}$ via a small regularization parameter $\varepsilon > 0$. In this paper, $\varepsilon$ is specified as $\left\| \mathbf{y} - \mathbf{g}(\mathbf{m}^{(n)}) \right\|_2^2$ to adaptively control the weighting matrix during the iterations. The basic idea behind this choice is to place more weight on the regularization term at early iterations of the Newton method where the iterates are far from the true solution, and to gradually decrease the effect of sparsity as the data misfit reduces at later iterations. We note that if $p$=2 in (P2), the enforced regularization is the classical $l_2$-norm constraint that minimizes the energy of the solution and does not promote sparsity; when $0 \le p \le 1$ is used, as shown in Fig. 2, the sparsity regularization is implicitly enforced where the $l_p$-quasinorm approximates the $l_0$-quasinorm in (P0).

When Newton method is used to minimize (P2), the updated values of $\mathbf{m}$ at iteration $n$+1 are obtained using the criterion

$$\mathbf{m}^{(n+1)} = \left\{ \mathbf{m} : \frac{\partial J_{n+1}(\mathbf{m})}{\partial \mathbf{m}} = \mathbf{0} \right\} \tag{9}$$

which leads to iterative solution of

$$\left( \mathbf{G}^T \mathbf{G} + \alpha \boldsymbol{\Phi}^T \mathbf{W} \boldsymbol{\Phi} \right) \mathbf{m}^{(n+1)} = \mathbf{G}^T \mathbf{y}_n \tag{10}$$

where $\mathbf{y}_n = \mathbf{y} - \mathbf{g}\left(\mathbf{m}^{(n)}\right) + \mathbf{G}\mathbf{m}^{(n)}$, and $G$ represents Jacobian matrix in the first-order Taylor approximation of $\mathbf{g}\left(\mathbf{m}\right)$ around $\mathbf{m}^{(n)}$, that is $\mathbf{g}\left(\mathbf{m}\right) \approx \mathbf{g}\left(\mathbf{m}^{(n)}\right) + \mathbf{G}\left(\mathbf{m} - \mathbf{m}^{(n)}\right)$. We obtain the Jacobian matrix $G$ by solving the adjoint model [REF] for each iteration. Unfortunately, choosing a suitable regularization parameter $\alpha$ is not straightforward. There are standard methods for specifying $\alpha$ in the literature such as the L-curve or the generalized cross validation (GCV) methods. However, these methods add a computational burden on the solution algorithm and are developed for linear inverse problems where the discrete Picard condition (REF) is satisfied for the underlying unperturbed problem. The additional computation can become prohibitive for nonlinear dynamic inverse problems with time-consuming forward simulation runs, especially when the initial guess is very far from the true solution. The summary of the algorithm used to solve (P2) is outlined in Table I.

### III. B Multiplicative Regularization

While adding the sparsity promoting penalty to the cost functional has a positive effect on the quality of the reconstructions, it also adds the unknown regularization parameter which is not trivial to specify. To eliminate the regularization parameter, we include the regularization term as a multiplicative constraint. In this approach, the importance of the regularization is *adaptively* determined during the inversion while removing the artificial regularization parameter. The resulting multiplicatively regularized objective function can be written as

$$J\left(\mathbf{m}\right) = \left\|\mathbf{y} - \mathbf{g}\left(\mathbf{m}\right)\right\|_2^2 \cdot SP\left(\mathbf{m}\right) \tag{11}$$

where $SP\left(\mathbf{m}\right) := \left\|\mathbf{\Phi}\mathbf{m}\right\|_W^2 = \mathbf{m}^T\mathbf{\Phi}^T\mathbf{W}\mathbf{\Phi}\mathbf{m}$ is the regularization term that promotes the sparsity of $\mathbf{\Phi}\mathbf{m}$. The diagonal "weighting" matrix $\mathbf{W}$ is similarly to the additive regularization case,

We note that the structure of the cost functional in (11) is such that it minimizes the penalty term with a large regularization weight at the early iterations of the Newton's algorithm, when the value of $\left\|\mathbf{y} - \mathbf{g}\left(\mathbf{m}\right)\right\|_2^2$ is large. As the iterations proceed and the data misfit term is reduced so is the weight on the sparsity regularization term. This adaptive behavior is amplified with the specified choice of the parameter $\varepsilon_n$, i.e. $\left\|\mathbf{y} - \mathbf{g}\left(\mathbf{m}^{(n)}\right)\right\|_2^2$; that is, during early iterations (larger observation misfit) the weighting matrix $\mathbf{W}_{i,i} = \left( \left( \mathbf{\Phi}\mathbf{m}^{(n)} \right)_i^2 + \varepsilon_n \right)^{\frac{p}{2}-1}$ gives a larger weight to the sparsity regularization (note that $0 \leq p \leq 1$) while as the iterations proceed (and the data misfit is reduced) less weight is adaptively placed on the regularization term. which eliminates the need for an additional regularization parameter. Based on the simulation results in the next section, the multiplicative regularization approach seems to perform well even when relatively limited noisy data are used.

The ($n$+1)-th iteration of the Netwon's method under multiplicative regularization is derived from Equation (9) as follows

$$-\mathbf{G}^T \left( \mathbf{y} - \mathbf{g}(\mathbf{m}) \right) SP(\mathbf{m}) + \left\| \mathbf{y} - \mathbf{g}(\mathbf{m}) \right\|_2^2 \mathbf{\Phi}^T \mathbf{W} \mathbf{\Phi} \mathbf{m} = 0 \qquad (12)$$

Which leads to the update equation

$$\left( \mathbf{G}^T \mathbf{G} + \beta(\mathbf{m}) \mathbf{\Phi}^T \mathbf{W} \mathbf{\Phi} \right) \mathbf{m} = \mathbf{G}^T \left( \mathbf{y} - \mathbf{g}(\mathbf{m}^{(n)}) + \mathbf{G}\mathbf{m}^{(n)} \right) \qquad (13)$$

where $\beta(\mathbf{m}) = \dfrac{\left\| \mathbf{y} - \mathbf{g}(\mathbf{m}) \right\|_2^2}{SP(\mathbf{m})}$, and $g(\mathbf{m}) \approx g(\mathbf{m}^{(n)}) + G(\mathbf{m} - \mathbf{m}^{(n)})$ has been explicitly assumed. The solution algorithm for the multiplicative regularization is given in Table II.

## IV. RESULTS OF NUMERICAL EXPERIMENTS

In this section, two sets of synthetic numerical example resembling waterflooding experiments in an oil reservoir have been carried out to evaluate the performance of our algorithm. In each example a 320 m × 320 m × 10 m synthetic reservoir is discretized into a two-dimensional 32 × 32 × 1 uniform grid block system with 10m × 10 m × 10 m uniform block sizes (see Table III). The simulations are performed using a two-phase (oil/water) immiscible in-house simulator. The waterflooding experiments are designed to inject one pore volume of water into the reservoir within a one year simulation time. During the one year simulation time, 30 observation intervals of approximately every 12 days are considered and the measurements. Two different well configurations are assumed to study the performance of proposed method under varying amount of available observations. We refer to the first and second configurations in Fig. 3 as Reservoir A and Reservoir B, respectively. Reservoir A (Fig. 3a) portrays a line drive injection using a horizontal well with 32 ports that uniformly inject (from the left end of the domain) one pore volume of water into the reservoir during the one year simulation. A similar horizontal well with 32 ports is placed at the right end of the domain to produce the oil and water that move toward the production well. The production ports are under a total production rate constraint to preserve mass balance; these ports produce an equal volume of fluid (oil and water) as the volume of the injected water into the reservoir Reservoir B ( Fig. 3b) includes four injection wells (shown with filled black squares in Fig. 3b) and six production wells (shown with empty black squares). The injectors uniformly inject a total of one pore volume of water into the reservoir from the left end of the domain during the one year simulation and the producers at the right end of the reservoir uniformly extract an equal one pore volume of fluids. The initial and boundary conditions are assumed to be known perfectly and are listed, along with other important input parameters, in Table III. For all numerical tests, the initial solution for the permeability is homogeneous with a permeability of 20 mD. We have used the DCT as the compression transform domain (basis $\Phi$ is the full rank $32^2$ dimensional DCT basis) in which a sparse solution is sought.

**Numerical Examples 1: Reservoir A**

We first consider the configuration in Reservoir A and discuss the following three choices with the Iteratively Reweighting solution method.

- $p=2$ which implies that sparsity is not promoted in the solution;
- $p=1$ which results in the basis pursuit formulation (REF) where the $l_1$-norm convex relaxation is used to approximates the $l_0$-quasinorm;
- $p=0$ which leads to the strong sparsity promoting regularization via approximate $l_0$-quasinorm constraint.

To illustrate the reconstruction performance for different $p$ values, the reconstructions for $p=1.5$ and $p=0.5$ are also considered for this example.

Figure 4a displays the true permeability and the corresponding water saturation solution profiles after 0 months, 2 months, 4 months, 6 months and 12 months during the simulation. The preferential water flow inside the channels and the resulting early water breakthrough at the production wells are apparent in these figures. Figure 4b shows permeability solutions at sample iterations of the proposed minimization scheme for $p=2$, where a minimum energy regularization is used, resulting in a smooth solution. The maximum number of minimization iterations was set to 30; however, in most cases after 15 iterations no major improvements in the objective function and the estimated parameters were observed. The saturation profiles corresponding to the final permeability solution are shown in Fig. 4b. From Fig.4b, it is obvious that the reconstruction still stays near the initial homogeneous permeability guess, which is far from the true permeability. Hence, the estimated saturation profiles fail to capture the behavior of the fluid displacements and therefore can not predict the significant amount of bypassed oil in the reservoir. In practice, such reconstruction errors can adversely affect future field development activities such injection/production control strategies and well placements. In order to show this, Fig. 6 gives the samples of comparsions between true pressure/saturation and reconstructed ones for different $p$.

Figure 4c shows the iteration solution samples of the same inverse problem for $p$=1 where the the $l_1$-norm convex relation of the original problem in (P0) is used. All other conditions remain identical to the previous example. The saturation profiles corresponding to final permeability solution are also shown in Fig. 4c. As seen from Fig. 4c, the reconstructed solution effectively captures the location and orientation of the high permeability channels in this case Indicating that the $l_1$-norm sparsity constraint (i.e. basis pursuit algorithm) is better able to identify the shape and continuity of the true channels, even though observations are only collected from the two ends of the reservoir. As a result the corresponding saturation profiles in Fig. 4c identify the preferential water flow displacement behavior and the resulting bypassed oil, which is quite important in optimizing the future oil production strategies from the reservoir.

Figure 4d shows the sample iteration solution of the same inverse problem for $p$=0, where an approximate strong sparsity regularization is promoted via the $l_0$-norm minimization. The reconstruction has obviously failed in this case mainly because of the non-convex property from $l_0$-norm. The results during the first few iterations (iterations 1 and 4) look relatively more reasonable because during these iterations, $\left\| \mathbf{y} - \mathbf{g}\left(\mathbf{m}^{(n)}\right) \right\|_2^2$ is relatively large; consequently, the weighting matrix $\mathbf{W}_{i,i} = \left( \left( \mathbf{\Phi}\mathbf{m}^{(n)} \right)_i^2 + \varepsilon_n \right)^{\frac{p}{2}-1}$ reduces the effect of the strong sparsity constraint. With the reduction in $\left\| \mathbf{y} - \mathbf{g}\left(\mathbf{m}^{(n)}\right) \right\|_2^2$ at later iterations, the $l_0$-quasinorm constraint becomes more important in the cost function. The nonconvexity of the $l_0$-norm and the lack of smoothness in its gradient can bring about unpredictable convergence issues. Figure 5 shows that the cost function is not well behaved during the iterations and a monotonically decreasing trend is not observed.

Samples iterations and the final permeability solution for $p$=1.5 and $p$=0.5 are shown in Fig. 4e and 4f, respectively. The results indicate that that the reconstruction with $p$=1.5 has similar behavior to the the case with $p$=2 and does not promote solution sparsity whereas the reconstruction with $p$=0.5 produces similar reconstruction results as in the case with $p$=1. These results reiterate the effectiveness of the case with $p$=1, i.e. $l_1$-norm regularization, (known as the basis pursuit algorithm) in promoting sparsity.

**Numerical Examples 2: Reservoir B**

In this subsection, we use Reservoir B with fewer observations to study the effect of observations. An important theoretical requirement for successful performance of the compressed sensing for linear problems is the availability of sufficient measurements [REF]. Intuitively, this requirement (or possibly a stronger version of it) is also needed for non-linear problems, which is the case that we are investigating in this paper. It is however, important to note that the production observations are available at 30 time intervals in our examples, which substantial increases the nature and number of available measurements. Reservoir B has a total of 10 observation locations (wells) as opposed to Reservoir A with 64 measurement locations. The true permeability is similar to the one used in previous examples (see Figs. 4a-f). Figure 7a displays the true permeability and the corresponding water saturation solution profiles after 0 months, 2 months, 4 months, 6 months and 12 months during the simulation. We only present the result for the basis pursuit reconstruction method (i.e. $p$=1). Figure 7b shows the sample permeability iterations and the saturation profiles corresponding to the final solution for this case, respectively. In addition, the comparison between true pressure/saturation and reconstructed ones have been provided in Fig. 8. It can be seen that the reconstruction are quite promising even in this case where limited spatial measurements are available. An important observation can be made about the results in Figure 7; in particular, the channel structure is correctly captured even when there is no any prior information besides the compressibility of permeability distribution.

As the final example, a more realistic example about the reconstruction of heterogeneous reservoir is investigated under the condition of Reservoir B. Figure 9a displays the true permeability and the corresponding water saturation solution profiles after 0 months, 2 months, 4 months, 6 months and 12 months during the simulation. We only present the result for $p$=1. Figure 9b shows the sample permeability iterations and the saturation profiles corresponding to the final solution for this case, respectively. In addition, the comparison between true pressure/saturation and reconstructed ones have been provided in Fig. 9c. From these figures it can be seen that the position , range and value of high/low permeability can be clearly identified though non-perfect reconstruction.

## V. CONCLUSIONS

In this paper, a novel estimation approach is introduced for solving ill-posed nonlinear dynamic inverse problems with unknown parameters that are approximately sparse in a certain transform domain such as Fourier or wavelet domains. The proposed method is inspired by recent developments in sparse reconstruction literature and the celebrated compressed sensing theory. By recognizing that spatially correlated geologic facies and their flow-related properties (such as rock permeability) often have sparse representations/approximation in an appropriate compression transform domain, we search for a solution that is sparse in the specified transform domain and satisfies the dynamic flow measurements. The solution approach is formulated by minimizing an objective function that minimizes the data misfit and a multiplicative regularization term that promotes solution sparsity. We have presented several regularization methods and their corresponding reconstruction results, which suggest that among the regularization term studied minimization of a multiplicative $l_1$-norm regularization term is quite effective. While for linear problems, this choice is related to compressed sensing and basis pursuit algorithms that use a convex relaxation of the objective function, for our nonlinear problem the objective function is in general nonconvex and we have used a gradient-based solution technique to find a local solution. Nonetheless, the $l_1$-norm regularization appears to perform well in promoting solution sparsity.

We have applied the proposed algorithm to numerical examples that are drawn from reservoir characterization where the estimation of geologic facies was considered. The results we have presented in this paper show that the proposed approach can provide satisfactory reconstruction for nonlinear dynamic inverse problems studied. Our results support the hypothesis that *$l_1$-norm based sparsity promoting regularization* in a sparse domain can be used to regularize the nonlinear inverse problems in geosciences applications. In our example, the *$l_1$-norm* sparsity regularization could successfully identify the significant transform basis vectors and adaptively tune their corresponding coefficients (dynamic range) using available observations and the reconstruction results in previous iterations. In summary, the results presented in this paper suggest that sparsity constraints via $l_1$-norm minimization provide a promising approach for regularization of ill-posed inverse problems with approximately sparse unknowns in an appropriate transform domain. Since the only prior information we have used is the sparsity of the solution in a transform domain, we anticipate that the quality of the solution can be further improved if additional prior information is included in the solution. The proposed algorithm in this paper can have far reaching implications in solving many real-world nonlinear imaging inverse problems, including those encountered in geophysical, electrical, and medical tomography.

Proceedings of ICASSP'08, submitted.

**Table I**

**Table I. The solution procedure for additive regularization in (P2)**

---

**Initialization**:

1. Initialization $\boldsymbol{m}$

2. Select a suitable sparse transformation $\boldsymbol{\Phi}$

3. Specify $p$ and $\alpha$

5. Compute the weight matrix $\boldsymbol{W}$

      **While** (stopping criterion not met)

          **DO**

               Update $\alpha$ (optional)

               Solve $\left(\mathbf{G}^T\mathbf{G} + \alpha\boldsymbol{\Phi}^T\mathbf{W}\boldsymbol{\Phi}\right)\mathbf{m}^{(n+1)} = \mathbf{G}^T\mathbf{y}_n$

               Update $\boldsymbol{W}$

               Update sensitivity matrix $\boldsymbol{G}$

               Compute the convergence criterion

          **END**

---

**Table II**

**Table II. The pseudo-code of proposed algorithm**

---

**Initialization**:

1. Initialization $\boldsymbol{m}$

2. Select a suitable sparse transformation $\boldsymbol{\Phi}$

3. Specify $p$

5. Compute the weight matrix $\boldsymbol{W}$

      **While** (stopping criterion not met)

        **DO**

            Solve $\left(\mathbf{G}^T\mathbf{G}+\beta\left(\mathbf{m}\right)\boldsymbol{\Phi}^T\mathbf{W}\boldsymbol{\Phi}\right)\mathbf{m}=\mathbf{G}^T\left(\mathbf{y}-\mathbf{g}\left(\mathbf{m}^{(n)}\right)+\mathbf{G}\mathbf{m}^{(n)}\right)$

            Update $\boldsymbol{W}$

            Update sensitivity matrix $\boldsymbol{G}$

            Compute the convergence criterion

        **END**

**Table III**

**Table III: General simulation information for Reservoir A and B**

| Parameter | Reservoir A | Reservoir B |
|---|---|---|
| **Simulation Parameters** | | |
| Phases | Two-phase (o/w) | Two-phase (o/w) |
| Simulation Time | 1 year | 1 year |
| Grid systems | 32 by 32 by 1 | 32 by 32 by 1 |
| Cell dimensions | 10 by 10 by 10 | 10 by 10 by 10 |
| Rock porosity | 0.20 | 0.20 |
| Initial oil saturation | 0.90 | 0.90 |
| Injection volume | 1PV | 1PV |
| Number of injectors | 32 | 4 |
| Number of producers | 32 | 6 |
| **Assimilation Information** | | |
| Observation intervals | 12 days | 12 days |
| Obs. at injection wells | Pressure | Pressure |
| Obs. at production wells | Pressure& saturation | Pressure& saturation |

**FIGURE CAPTIONS**

**Figure 1.** The distribution of classical reservoir permeability, its DCT coefficients, and the reconstruction with 5% largest DCT coefficients.

**Figure 2.** The comparison of $|x|^p$ for different $p$.

**Figure 3.** The sketch map of two considered well configuration

**Figure 4.**

Figure 4(a) The true permeability and the saturation after 0months, 2 months, 4 months, 6 months and 12 months.

Figure 4(b-f) The reconstructed permeability of reservoir A with different $p$ after 1 iteration, 3 iterations, 6 iterations, 9 iterations, 12 iterations and 15 iterations and the corresponding saturation after 0 months, 2 months, 4 months, 6 months and 12 months.

The initial solution is homogeneous distribution with the permeability of 20 mD.

**Figure 5.** Convergence curve of figure.4d, where the $x$-axis is the iteration steps while $y$-axis is the error described by cost functional.

**Figure 6.** The comparison between the true and reconstructed pressure/saturation at production wells.

**Figure 7.**

Figure 7(a) The true permeability and the saturation after 0months, 2 months, 4 months, 6 months and 12 months.

Figure 7(b) The reconstructed permeability of reservoir A with $p=1$ after 1 iteration, 3 iterations, 6 iterations, 9 iterations, 12 iterations and 15 iterations and the corresponding saturation after 0 months, 2 months, 4 months, 6 months and 12 months.

The initial solution is homogeneous distribution with the permeability of 20 mD.

**Figure 8.** The comparison between the true and reconstructed pressure/saturation at production wells.

**Figure 9.**

Figure 9(a) The true permeability and the saturation after 0months, 2 months, 4 months, 6 months and 12 months.

Figure 9(b) The reconstructed permeability of reservoir A with $p=1$ after 1 iteration, 3 iterations, 6 iterations, 9 iterations, 12 iterations and 15 iterations and the corresponding saturation after 0

months, 2 months, 4 months, 6 months and 12 months.

The initial solution is homogeneous distribution with the permeability of 20 mD.

Figure 9(c) The comparison between the true and reconstructed pressure/saturation at production wells.

**FIGURES**

## Figure 1

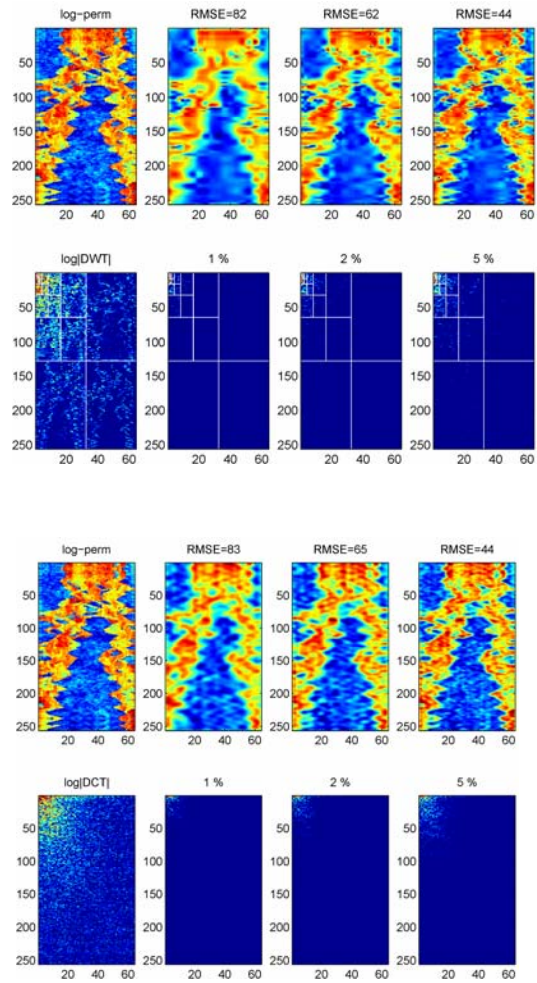

**Figure 2**

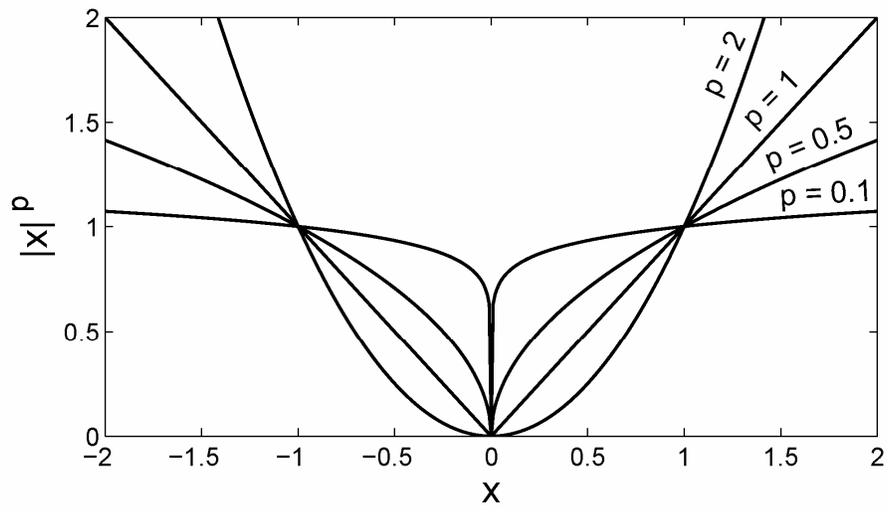

**Figure 3**

(a)                                      (b)

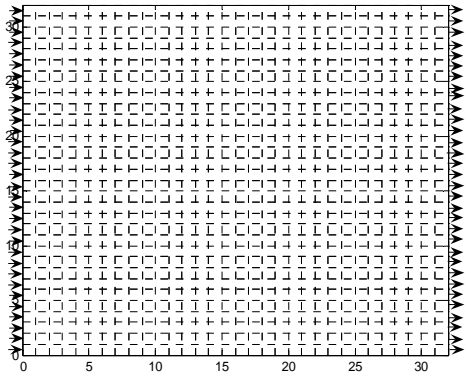          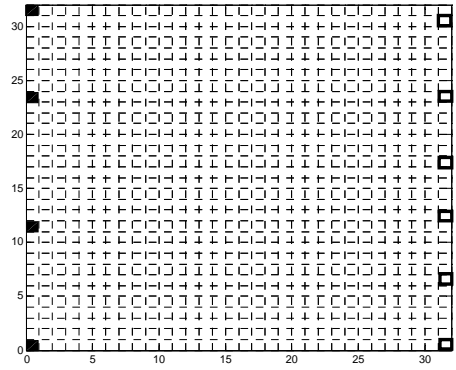

# Figure 4

## Fig. 4(a) True perm and corresponding saturation profiles

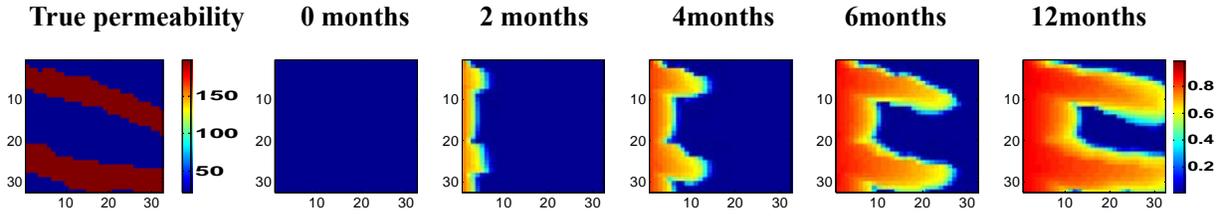

## Fig. 4(b) Estimated perm and corresponding saturation profiles *(p=2)*

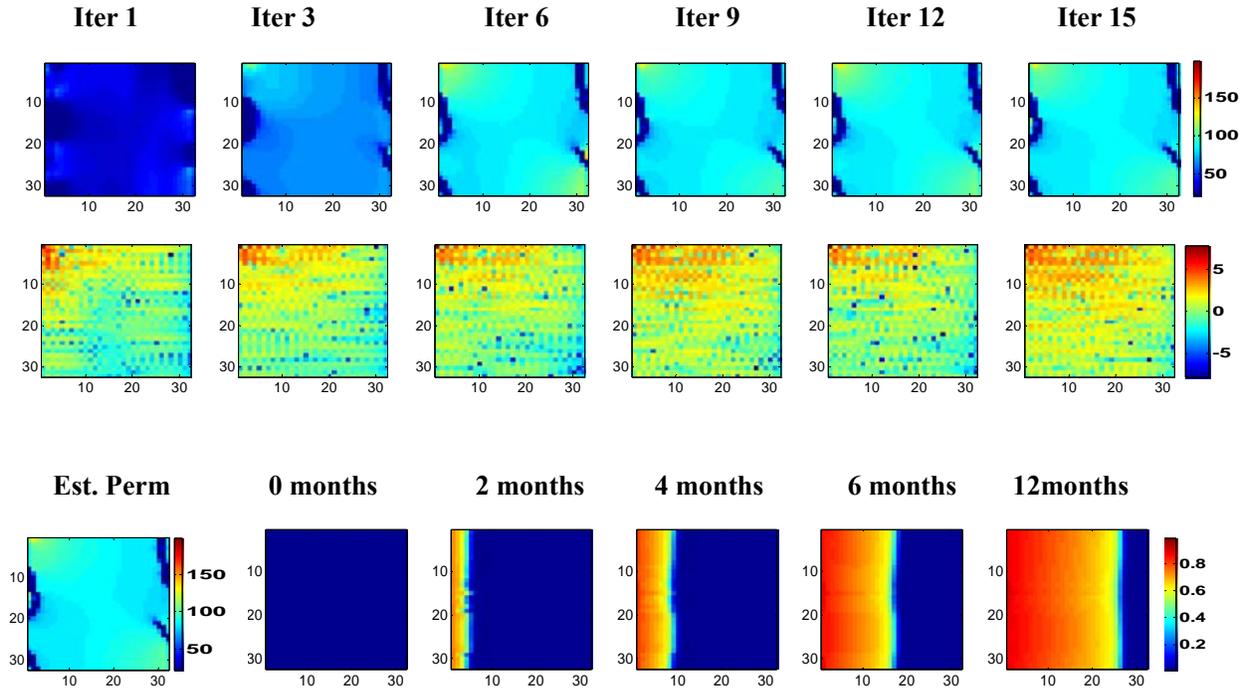

## Fig. 4(c) Estimated perm and corresponding saturation profiles *(p=1)*

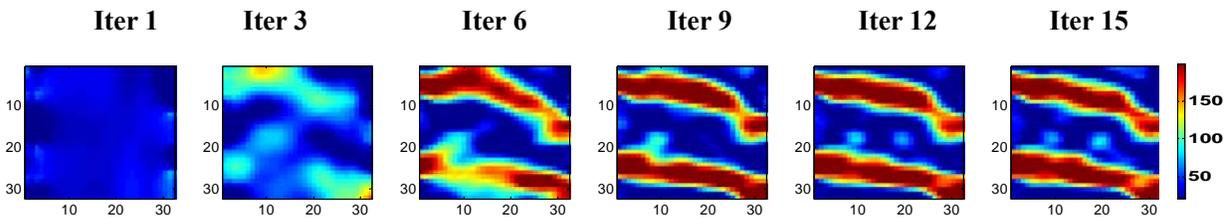

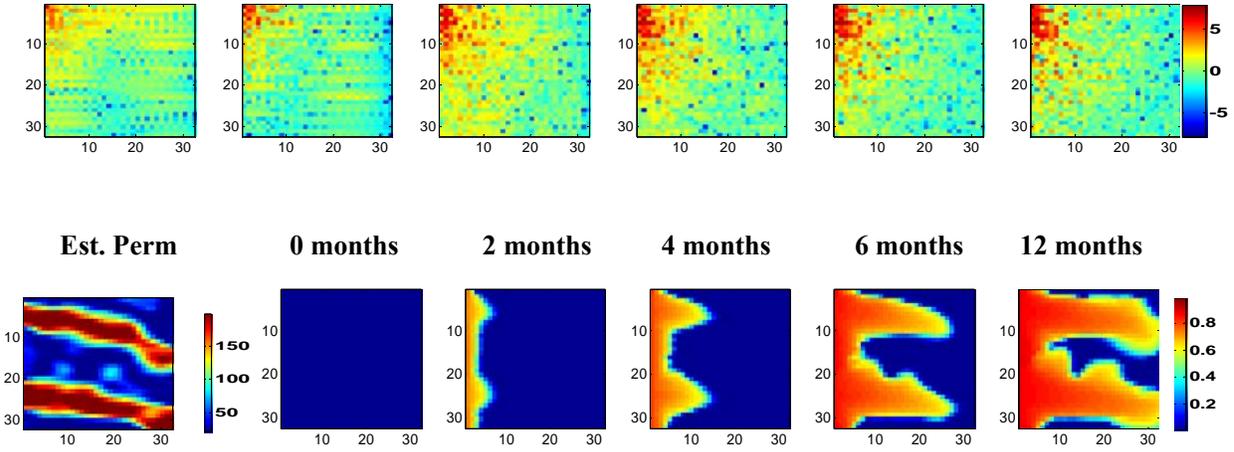

**Est. Perm**   **0 months**   **2 months**   **4 months**   **6 months**   **12 months**

**Fig. 4(d)** Estimated perm and corresponding saturation profiles *(p=0)*

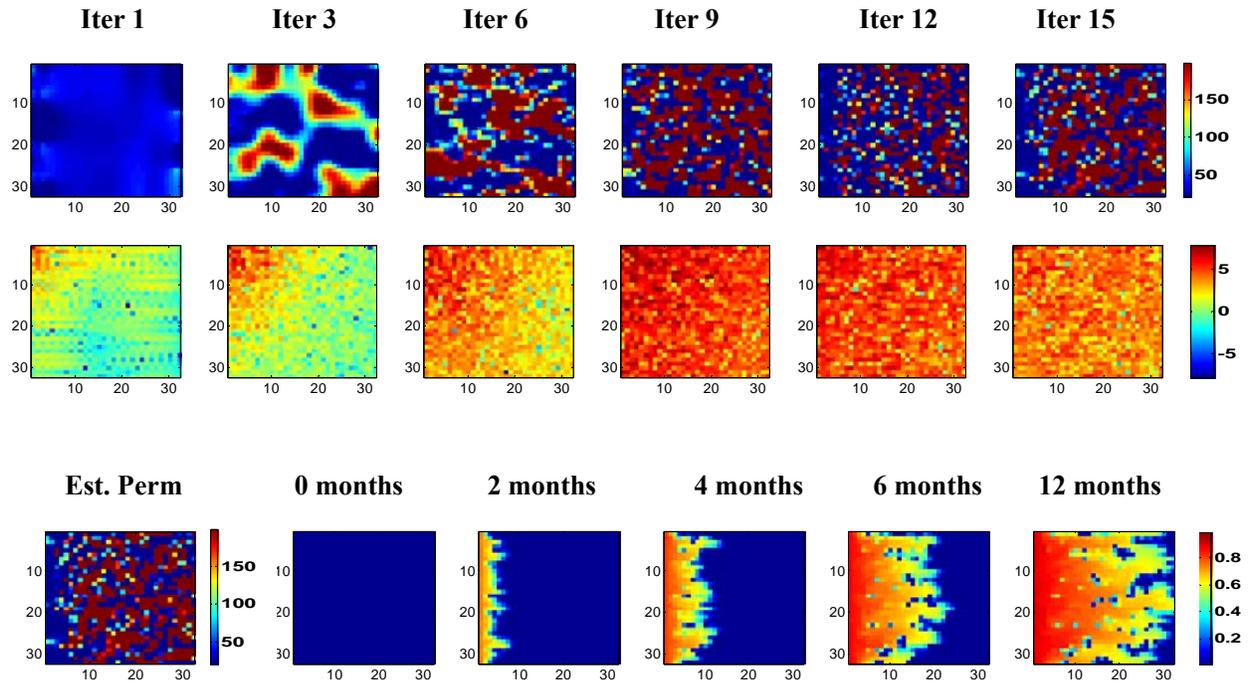

**Iter 1**   **Iter 3**   **Iter 6**   **Iter 9**   **Iter 12**   **Iter 15**

**Est. Perm**   **0 months**   **2 months**   **4 months**   **6 months**   **12 months**

**Fig. 4(e)** Estimated perm *(p=1.5)*

**Iter 1**   **Iter 3**   **Iter 6**   **Iter 9**   **Iter 12**   **Iter 15**

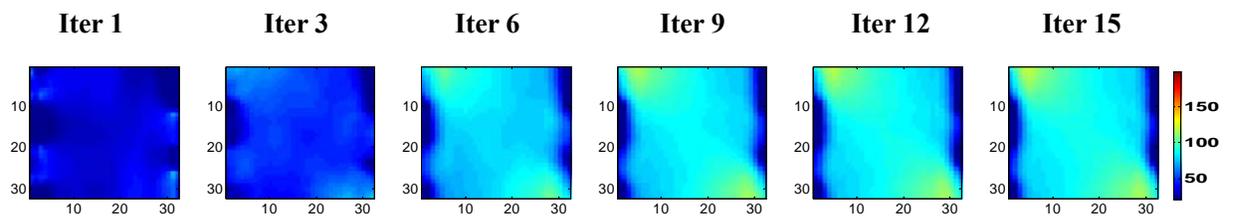

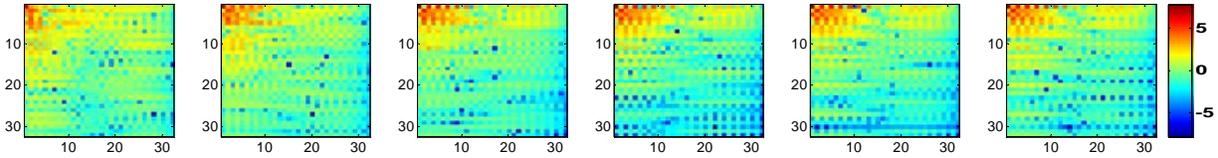

**Fig. 4(f)** Estimated perm *(p=0.5)*

|  Iter 1 | Iter 3 | Iter 6 | Iter 9 | Iter 12 | Iter 15 |

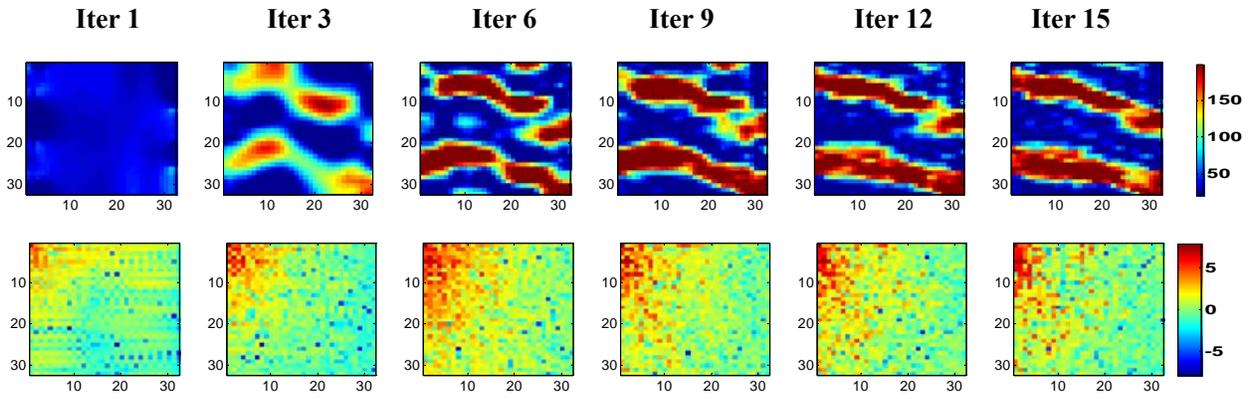

## Figure 5

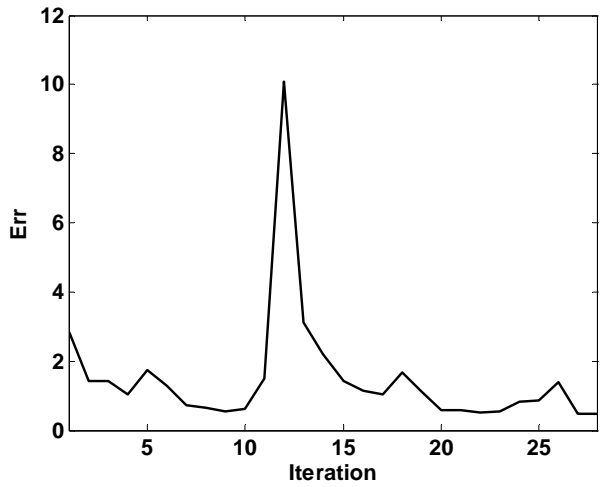

**Figure 6**

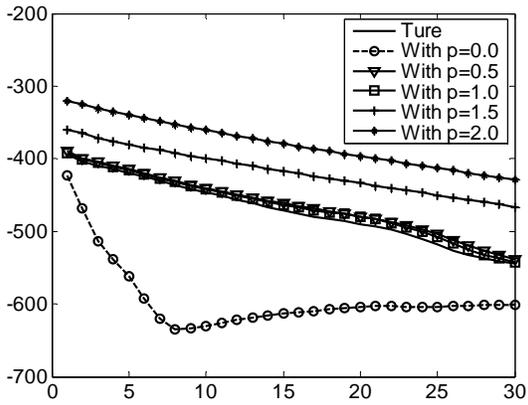

**(a) Pressure at prodcution well #1**

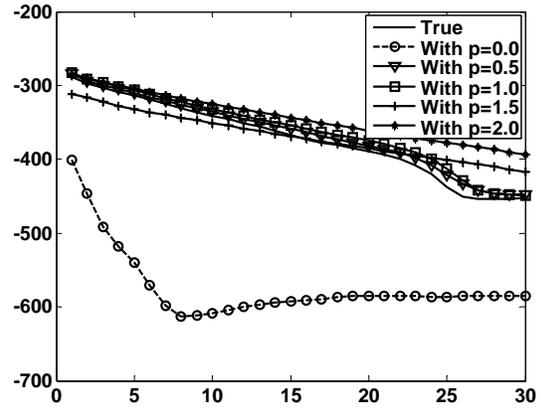

**(b) Pressure at production well # 10**

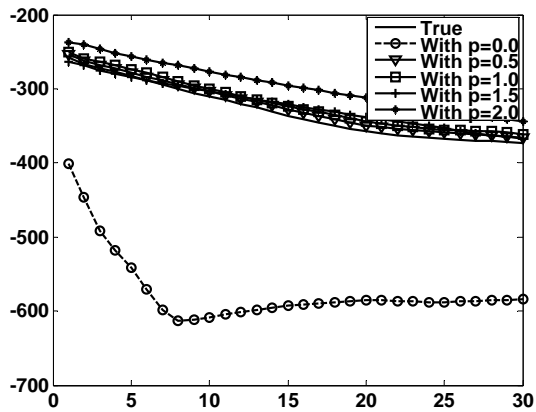

**(c) Pressure at prodcution well #20**

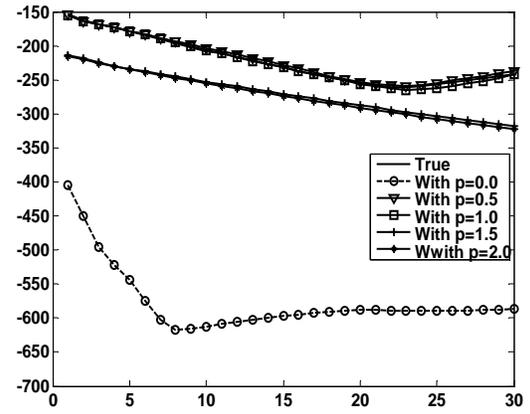

**(d) Pressure at production well # 30**

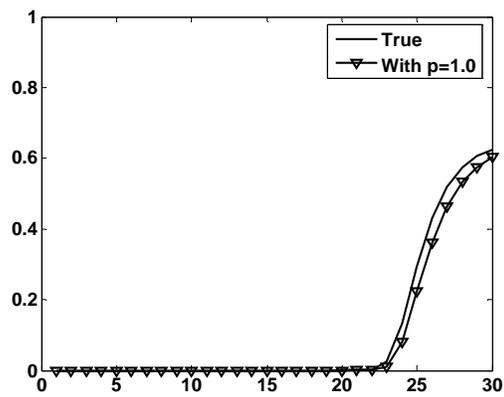

**(e) Satruation at production well #10,**

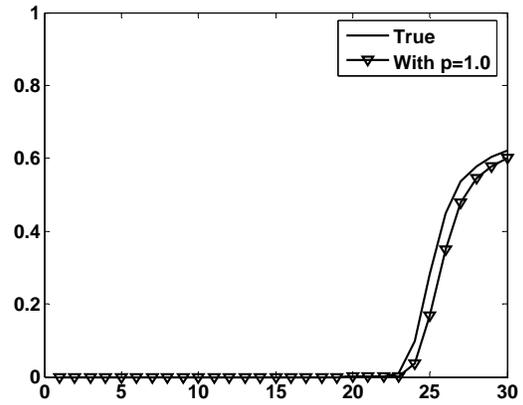

**(f) Saturation at production well #20**

# Figure 7

**Fig. 7a** True perm and corresponding saturation profiles

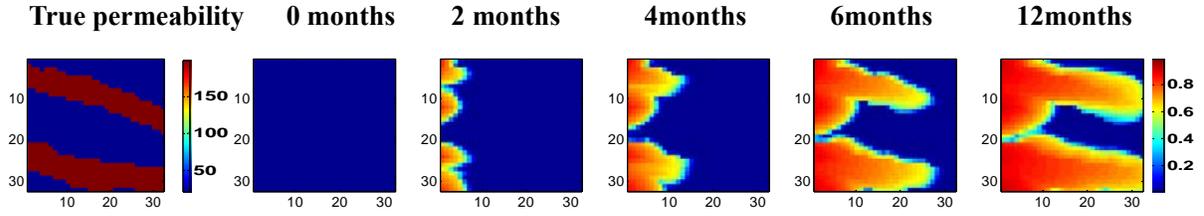

**Fig. 7(b)** Estimated perm and corresponding saturation profiles *(p=1)*

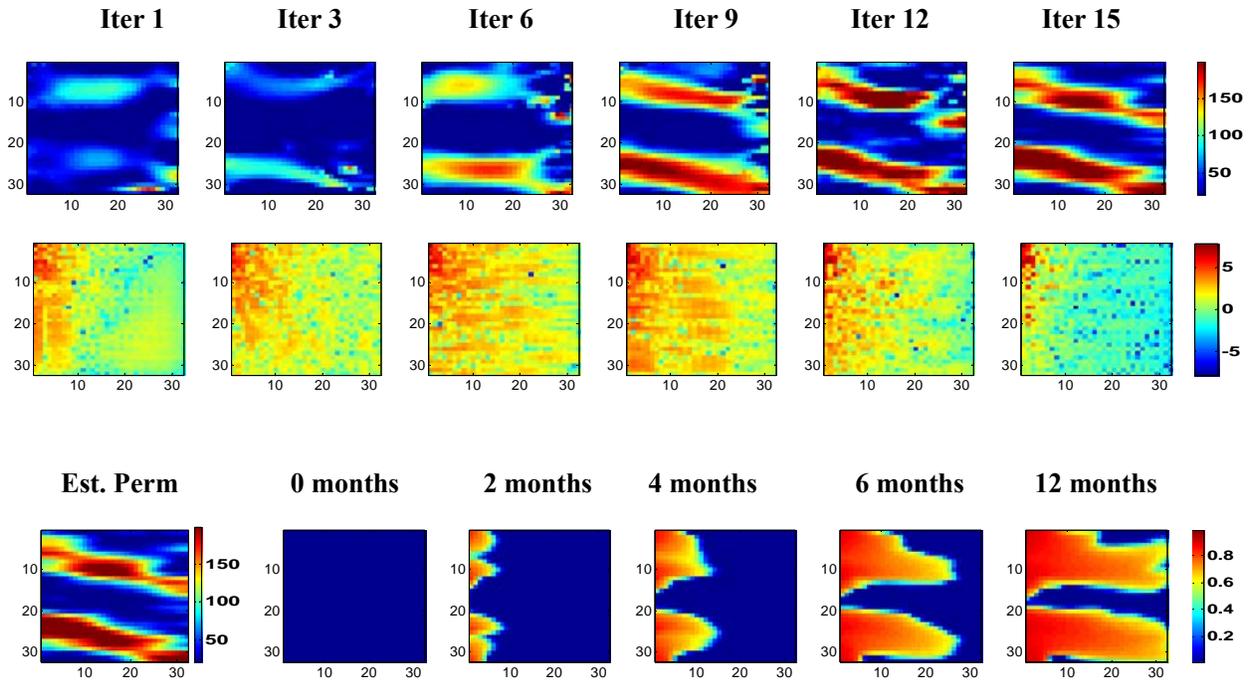

**Figure 8**

**Fig. 8a    True pressure(solid line) and reconstructed pressure (dashed line) at production wells**

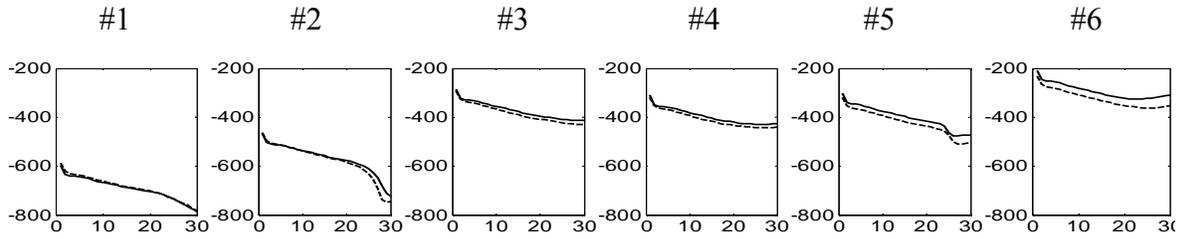

**Fig. 8b    True saturation and reconstructed saturation at some production wells**

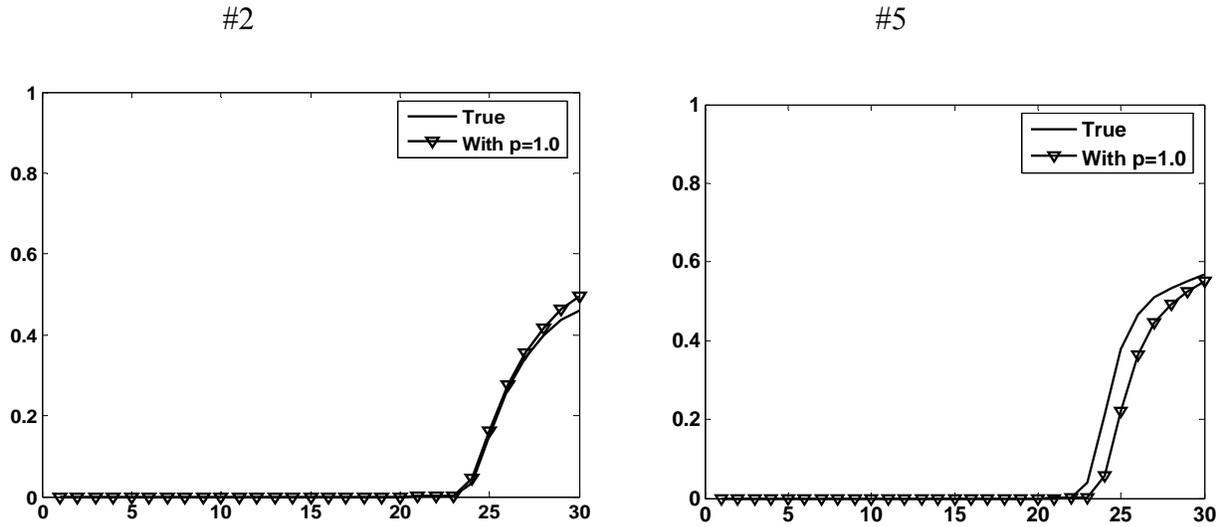

# Figure 9

### Fig. 9a True log-perm and corresponding saturation profiles

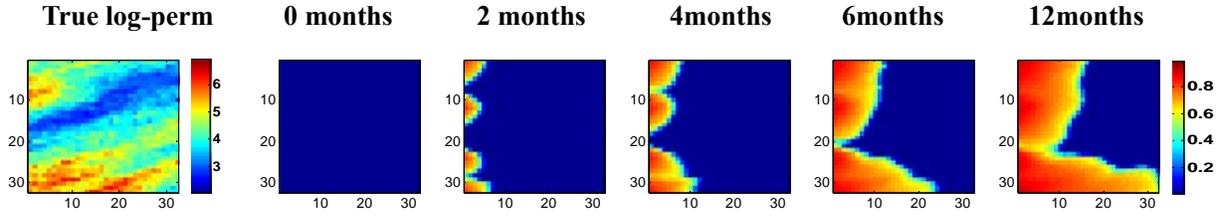

### Fig. 9(b) Estimated perm and corresponding saturation profiles *(p=1)*

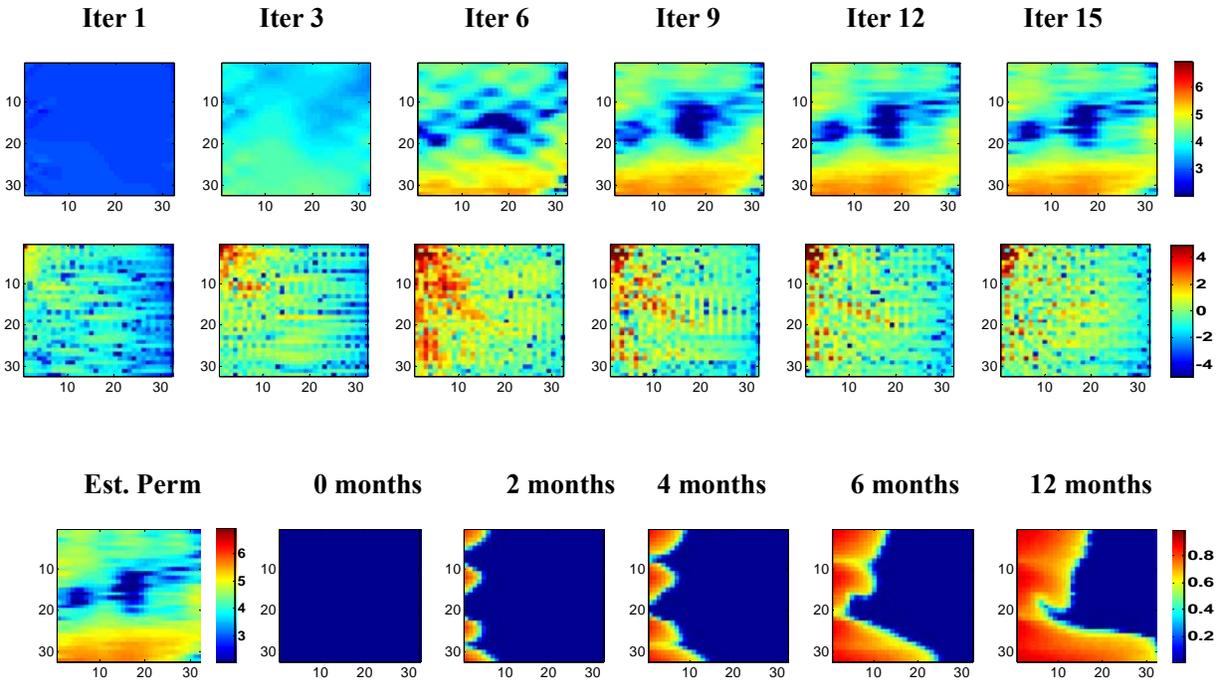

### Fig. 9(c) True pressure(solid line) and reconstructed pressure (dashed line) at production wells

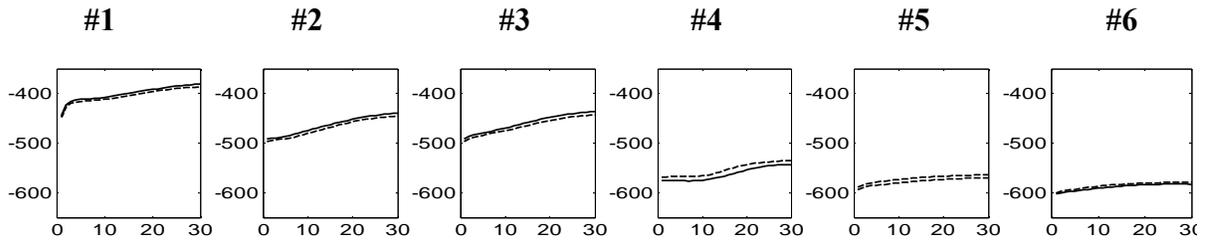